\documentclass[fleqn]{mat01}
\usepackage{times,mathtimy,amssymb,boxit}
\begin{document}

\setcounter{page}{91}
\firstpage{91}

\font\xx=msam5 at 9pt
\def\ab{\mbox{\xx{\char'03}}}

\font\sa=tibi at 10.4pt
\def\d{\hbox{d}}

\def\thoe{\trivlist\item[\hskip\labelsep{{\bf Theorem}}]}
\newtheorem{theo}{\bf Theorem}
\renewcommand\thetheo{\arabic{theo}}
\newtheorem{theor}[theo]{\bf Theorem}
\newtheorem{definit}[theo]{\rm DEFINITION}
\newtheorem{lem}[theo]{Lemma}

\newcommand{\R}{\mbox{$\mathbb{R}$}}
\newcommand{\N}{\mbox{$\mathbb{N}$}}
\newcommand{\Q}{\mbox{$\mathbb{Q}$}}
\newcommand{\C}{\mbox{$\mathbb{C}$}}
\newcommand{\Z}{\mbox{$\mathbb{Z}$}}
\newcommand{\K}{\mbox{$\mathbb{K}$}}

\title{Some functional equations originating from number theory}

\markboth{Soon-Mo Jung and Jae-Hyeong Bae}{Functional equations}

\author{SOON-MO JUNG and JAE-HYEONG BAE$^{*}$}

\address{Mathematics Section, College of Science and Technology, Hong-Ik University, 339-701 Chochiwon, Korea\\
\noindent $^*$Department of Mathematics, Chungnam National University, 305-764 Daejon, Korea\\
\noindent E-mail: smjung@wow.hongik.ac.kr; jhbae@math.cnu.ac.kr}

\volume{113}

\mon{May}

\parts{2}

\Date{MS received 15 September 2002}

\begin{abstract}
We will introduce new functional equations (\ref{eq:hana}) and
(\ref{eq:hana2}) which are strongly related to well-known formulae
(\ref{eq:young}) and (\ref{eq:young2}) of number theory, and investigate
the solutions of the equations. Moreover, we will also study some
stability problems of those equations.
\end{abstract}


\maketitle

\section{Introduction}

In 1940, Ulam gave a wide ranging talk before the Mathematics Club
of the University of Wisconsin in which he discussed a number of
important unsolved problems \cite{14}. Among those was the
question concerning the stability of homomorphisms:\vspace{.5pc}

Let $G_{1}$ be a group and let $G_{2}$ be a metric group with a metric
\hbox{$d(\cdot\,\,,\cdot)$}. Given any $\varepsilon > 0$, does there exist a
$\delta > 0$ such that if a function $h : G_{1} \to G_{2}$ satisfies the
inequality $d(h(xy),h(x)h(y)) < \delta$ for all $x,y\in G_{1}$, then
there exists a homomorphism $H : G_{1} \to G_{2}$ with $d(h(x),H(x)) <
\varepsilon$ for all $x \in G_{1}$?\vspace{.5pc}

If the answer is affirmative, the functional equation for homomorphisms
is said to be stable in the sense of Hyers and Ulam because the first
result concerning the stability of functional equations was presented by
Hyers. Indeed, he has answered the question of Ulam for the case
where $G_{1}$ and $G_{2}$ are assumed to be Banach spaces (see \cite{4}).

We may find a number of papers concerning the stability results of
various functional equations (see \cite{b,blz,1,2,3,g,gs,5,6,7,8,9,10,11,12}
and the references cited therein).

According to a well-known theorem in number theory, a positive integer
of the form $m^2 n$, where each divisor of $n$ is not a square of the
integer, can be represented as a sum of two squares of integer if and
only if every prime factor of $n$ is not of the form $4k + 3$. In the
proof of this theorem, we make use of the following elementary
equalities
\begin{equation}
(x_1^2 + y_1^2)(x_2^2 + y_2^2) =
(x_1 x_2 + y_1 y_2)^2 + (x_1 y_2 - y_1 x_2)^2
\label{eq:young}
\end{equation}
and
\begin{align}
&(x_1^2 + y_1^2 + z_1^2 + w_1^2)(x_2^2 + y_2^2 + z_2^2 +
w_2^2) =  (x_1 x_2 + y_1 y_2 + z_1 z_2 + w_1 w_2)^2\nonumber\\
&\quad \ +\, (x_1 y_2 - y_1 x_2 + z_1 w_2 - w_1 z_2)^2
  +\, (x_1 z_2 - y_1 w_2 - z_1 x_2 + w_1 y_2)^2 \nonumber\\
&\quad \ +\, (x_1 w_2 + y_1 z_2 - z_1 y_2 - w_1 x_2)^2.\label{eq:young2}
\end{align}
As we know, the above equations explain that the product of any sums of
two (four) squares of integer is also a sum of two (four) squares of
integer.

These equalities (\ref{eq:young}) and (\ref{eq:young2}) may be
formulated by the following functional equations
\begin{equation}
f(x_1, y_1) f(x_2, y_2) = f(x_1 x_2 + y_1 y_2, x_1 y_2 - y_1 x_2)
\label{eq:hana}
\end{equation}
and
\begin{align}
&f(x_1, y_1, z_1, w_1) f(x_2, y_2, z_2, w_2)
= f(x_1 x_2 + y_1 y_2 + z_1 z_2 + w_1 w_2,\,\nonumber\\
 &\quad \              x_1 y_2 - y_1 x_2 + z_1 w_2 - w_1 z_2,
      x_1 z_2 - y_1 w_2 - z_1 x_2 + w_1 y_2,\,\nonumber\\
&\quad \               x_1 w_2 + y_1 z_2 - z_1 y_2 - w_1 x_2).
\label{eq:hana2}
\end{align}

In this paper, the solutions and stability problems of the above
equations will be investigated.

\section{Solutions and stability of (\ref{eq:hana})}

We will first investigate the solutions of the functional equation
(\ref{eq:hana}) in the class of functions $f : \R^2 \to \R$.

\begin{theo}[\!]
If a function $f : \R^2 \to \R$ satisfies the
functional equation $(\ref{eq:hana})$ for all $x_1, x_2, y_1, y_2 \in
\R${\rm ,} then there exist a multiplicative function $m : \R \to \R$ and a
signum function $\sigma : \R^2 \to \{ \pm 1 \}$ such that
\begin{equation*}
f(x, y) =
\sigma(x, y)\, m \left( \sqrt{x^2 + y^2\,}\,\right)
\end{equation*}
 for all real numbers $x$ and $y$.
\end{theo}

\begin{proof}
Put $y_1 = y_2 = 0$ in (\ref{eq:hana}) to get
\begin{equation}
f(x_1, 0) f(x_2, 0) = f(x_1 x_2, 0)
\label{eq:dool}
\end{equation}
for all $x_1, x_2 \in \R$. Replace the $x_i$s by $x$ and the $y_i$s by
$y$ in (\ref{eq:hana}) to get
\begin{equation}
f(x, y) f(x, y) = f(x^2 + y^2, 0),
\label{eq:set}
\end{equation}
for any $x, y \in \R$. Using (\ref{eq:set}) twice, we have
\begin{equation*}
f(y, x) f(y, x) = f(y^2 + x^2, 0) = f(x^2 + y^2, 0) = f(x, y) f(x, y)
\end{equation*}
and hence we may define a function $\sigma_1 : \R^2 \to \{ \pm 1 \}$ by
\begin{equation}
f(y, x) = \sigma_1(x, y) f(x, y)
\label{eq:net}
\end{equation}
for all real numbers $x$ and $y$.

From (\ref{eq:hana}), (\ref{eq:set}) and (\ref{eq:net}), it follows that
\begin{align}
f(2xy, x^2 - y^2)
& =  f(x, y) f(y, x)\nonumber\\
& =  \sigma_1(x, y) f(x, y) f(x, y)
 =  \sigma_1(x, y) f(x^2 + y^2, 0)\label{eq:yosot}
\end{align}
for any real numbers $x$ and $y$.\pagebreak

We notice that for any given $u, v \in \R$, the following system
of equations
\begin{equation}
\begin{cases}
          2xy = u,\\[.3pc]
          x^2 - y^2 = v
\end{cases}
\label{eq:ilgob}
\end{equation}
has solutions $(x(u,v), y(u,v))$ in $\R^2$ as we see in the
following:
\begin{equation*} 
\hskip -2pc(x(u,v), y(u,v)) =
\begin{cases}
( \pm\sqrt{v}, 0 )    &\mbox{for~$u = 0$~and~$v \geq 0$},\\[.5pc]
( 0, \pm\sqrt{-v} \,) &\mbox{for~$u = 0$~and~$v < 0$},\\[.5pc]
\left( \pm\sqrt{\frac{v+\sqrt{u^2 + v^2}}{2}},
\pm\frac{u}{\sqrt{2v+2\sqrt{u^2 + v^2}}}\right)  &\mbox{for~$u \neq 0$.}
\end{cases}
\end{equation*}
It follows from (\ref{eq:ilgob}) that $x^2 + y^2 = \displaystyle\sqrt{u^{2} + v^{2}}$.
According to (\ref{eq:yosot}) and (\ref{eq:ilgob}), we obtain
\begin{equation}
f(u, v) = \sigma_1(x(u,v), y(u,v)) f\left( \sqrt{u^2 + v^2}, 0 \right)
\label{eq:yodolb}
\end{equation}
for any $u, v \in \R$.

Taking (\ref{eq:yodolb}) into account, we may introduce another function
$\sigma : \R^2 \to \{ \pm 1 \}$ that satisfies the equality
\begin{equation}
f(u, v) = \sigma(u, v) f\left( \sqrt{u^2 + v^2}, 0 \right)
\label{eq:ahob}
\end{equation}
for all $u, v \in \R$.

Finally, define a function $m : \R \to \R$ by $m(x) = f(x, 0)$ for each
$x \in \R$. Then, (\ref{eq:dool}) and (\ref{eq:ahob}) ensure that $m$ is
a multiplicative function and that
\begin{equation*}
f(x, y) = \sigma(x, y) \, m \Big(\sqrt{x^2 + y^2\,}\,\Big)
\end{equation*}
for all real numbers $x$ and $y$.
\end{proof}

We will now investigate some stability problem of the functional
equation (\ref{eq:hana}). In view of Theorem 1, we can guess that the
stability of (\ref{eq:hana}) is strongly connected with multiplicative
functions.

\begin{theo}[\!]
Let $X$ be a field and $M_1, M_2, N_1, N_2 : X \to
[0, \infty)$ be functions. If a function $f : X^2 \to \C$ satisfies the
following inequality
\begin{align}
&   | f(x_1, y_1) f(x_2, y_2) - f(x_1 x_2 + y_1 y_2, x_1 y_2 - y_1 x_2)
|\nonumber\\[.5pc]
&\quad \ \leq \min\{ M_1(x_1), M_2(x_2), N_1(y_1), N_2(y_2) \}
\label{eq:giyog}
\end{align}
for all $x_1, x_2, y_1, y_2 \in X${\rm ,} then $f(x,0)$ is either bounded or
multiplicative and further it satisfies
\begin{equation*}
| f(x,y)^2 - f( x^2 + y^2, 0 ) | \leq
  \min\{ M_1(x), M_2(x), N_1(y), N_2(y) \}
\end{equation*}
for any $x, y \in X$.
\end{theo}

\begin{proof}
With $y_1 = y_2 = 0$, (\ref{eq:giyog}) implies

\begin{equation*}
| f(x_1, 0) f(x_2, 0) - f(x_1 x_2, 0) | \leq
  \min\{ M_1(x_1), M_2(x_2), N_1(0), N_2(0) \}
\end{equation*}
for $x_1, x_2 \in X$. If we substitute $m(x)$ instead of $f(x, 0)$ in the above
inequality, then we have
\begin{equation*}
| m(x_1) m(x_2) - m(x_1 x_2) | \leq
  \min\{ M_1(x_1), M_2(x_2), N_1(0), N_2(0) \}
\end{equation*}
for all $x_1, x_2 \in X$.

Applying a theorem of Sz\'{e}kelyhidi \cite{13} (see Corollary 8.4
in \cite{8}), we conclude that $m$ is either bounded or multiplicative.

Finally, put $x_1 = x_2 = x$ and $y_1 = y_2 = y$ in (\ref{eq:giyog})
to get
\begin{equation*}
| f(x,y)^2 - f( x^2 + y^2, 0 ) | \leq
  \min\{ M_1 (x), M_2 (x), N_1 (y), N_2 (y) \}
\end{equation*}
for all $x, y \in X$.
\end{proof}

\section{Solutions and stability of (\ref{eq:hana2})}

We first prove a lemma which turns out to be indispensable for the
investigation of solutions of the functional equation (\ref{eq:hana2}).

\begin{lem}
For any given $a, b, c, d \in \R,$ the system of
equations
\begin{equation*}
\begin{cases}
            (x+z)(y+w) = a,\vspace{1mm}\\
            2xz-y^2-w^2 = b,\vspace{1mm}\\
            (x+z)(w-y) = c,\vspace{1mm}\\
            x^2 - z^2 = d
\end{cases}
\end{equation*}
has at least one solution $(x, y, z, w)$ in $\R^4$.
\end{lem}

\begin{proof}$\left.\right.$\vspace{.5pc}

\begin{enumerate}
\renewcommand{\labelenumi}{(\alph{enumi})}
\item If $a = c = d = 0$ and $b \leq 0$, then $(x, y, z, w)
= \left(0, \displaystyle\sqrt{-b/2}, 0, \displaystyle\sqrt{-b/2}\right)$ is a solution of our system of
equations.

\item If $b > 0$ and $d = 0$, set $x = z = \alpha \neq 0$ and we will
determine the value of $\alpha$ later. It follows from the first and
third equations that
\begin{equation*}
y = \frac{a-c}{4\alpha} ~~\mbox{and}~~ w = \frac{a+c}{4\alpha}.
\end{equation*}
By the second one, we get a biquadratic equation
\begin{equation*}
16\alpha^4 - 8b\alpha^2 - a^2 - c^2 = 0,
\end{equation*}
and one of its solutions is
\begin{equation*}
\alpha = \frac{\sqrt{b + \sqrt{a^2 + b^2 + c^2}}}{2} > 0.
\end{equation*}
Hence, the system of equations is solvable in $\R^4$ when $b > 0$ and $d
= 0$.

\item For the remaining cases under the condition $d = 0$: either if $a
= 0$, $b \leq 0$, $c \neq 0$ and $d = 0$, or if $a \neq 0$, $b \leq 0$,
$c = 0$ and $d = 0$, or if $a \neq 0$, $b \leq 0$, $c \neq 0$ and $d =
0$, then we follow the lines in part (b) and find out one solution of
our system of equations.

\item If $d \neq 0$, set $x = \sqrt{d+\alpha}$ and $z = \sqrt{\alpha}$
for some $\alpha \geq \max\{ 0, -d \}$ ($\alpha$ will be determined
later). By the first and third equations, we have
\begin{equation*}
y = \frac{a-c}{2\left(\sqrt{d+\alpha} + \sqrt{\alpha}\right)}
  ~~\mbox{and}~~
  w = \frac{a+c}{2\left(\sqrt{d+\alpha} + \sqrt{\alpha}\right)}.
\end{equation*}
If we substitute those expressions for $x, y, z, w$ in the second one
and if we carry out a tedious calculation, then we get a quadratic
equation
\begin{align*}
q(\alpha) &= 16(a^2 + c^2 + d^2)\alpha^2 +
              8(2d(a^2 + c^2 + d^2) - b(a^2 + c^2))\alpha\\[.5pc]
&\quad \ -
              (a^2 + c^2 + 2bd)^2 = 0.
\end{align*}
This equation has one solution $\alpha$ which is not less than $0$ and
$-d$ because of $q(0) \leq 0$ and $q(-d) = -(a^2 + c^2 - 2bd)^2 \leq 0$.
Thus, the system is solvable in $\R^4$ for $d \neq 0$. \vspace{-.5cm}
\end{enumerate}
\end{proof}

In the following theorem, we investigate the solutions of the
functional equation (\ref{eq:hana2}) by the same idea that was
applied to the proof of Theorem~1.

\begin{theo}[\!]
If a function $f : \R^4 \to \R$ satisfies the
functional equation $(\ref{eq:hana2})$ for all $x_i, y_i, z_i, w_i \in
\R$ $(i = 1, 2),$ then there exist a multiplicative function $m : \R \to
\R$ and a signum function $\sigma : \R^4 \to \{ \pm 1 \}$ such that
\begin{equation*}
f(x, y, z, w) = \sigma(x, y, z, w)\, m \left( \sqrt{x^2 + y^2 + z^2 +
w^2\,}\,\right)
\end{equation*}
for all real numbers $x, y, z, w$.
\end{theo}

\begin{proof}
If we set $y_i = z_i = w_i = 0$ ($i = 1, 2$) in
(\ref{eq:hana2}), then
\begin{equation}
f(x_1, 0, 0, 0) f(x_2, 0, 0, 0) = f(x_1 x_2, 0, 0, 0)
\label{eq:dool2}
\end{equation}
for all $x_1, x_2 \in \R$.
If we substitute $x, y, z, w$ for the $x_i, y_i, z_i, w_i$ in
(\ref{eq:hana2}), then we have
\begin{equation}
f(x, y, z, w) f(x, y, z, w) = f(x^2 + y^2 + z^2 + w^2, 0, 0, 0)
\label{eq:set2}
\end{equation}
for any $x, y, z, w \in \R$. Use eq.~(\ref{eq:set2}) twice to get
\begin{align*}
f(y, z, w, x) f(y, z, w, x)
& = f(y^2 + z^2 + w^2 + x^2, 0, 0, 0)\\
& = f(x^2 + y^2 + z^2 + w^2, 0, 0 ,0)\\
& = f(x, y, z, w) f(x, y, z, w).
\end{align*}
Therefore, we may define a function $\sigma_1 : \R^4 \to \{ \pm 1 \}$
by
\begin{equation}
f(y, z, w, x) = \sigma_1(x, y, z, w) f(x, y, z, w)
\label{eq:net2}
\end{equation}
for all $x, y, z, w \in \R$.

It follows from (\ref{eq:hana2}), (\ref{eq:set2}) and (\ref{eq:net2})
that
\begin{align}
&f((x+z)(y+w), 2xz-y^2-w^2, (x+z)(w-y), x^2 - z^2)\nonumber\\[.7pc]
&\quad \ = f(x, y, z, w) f(y, z, w, x)\nonumber\\[.7pc]
&\quad \ = \sigma_1(x, y, z, w) f(x, y, z, w) f(x, y, z, w)\nonumber\\[.7pc]
&\quad \ = \sigma_1(x, y, z, w) f(x^2 + y^2 + z^2 + w^2, 0, 0, 0)
\label{eq:yosot2}
\end{align}
for any real numbers $x, y, z, w$.

According to Lemma 3, we can easily see that
\begin{align*}
&\{ ((x+z)(y+w), 2xz-y^2-w^2, (x+z)(w-y), x^2 - z^2) \,:\,\\
&\quad x, y, z, w \in \R \} = \R^4
\end{align*}
because the following system of equations
\begin{equation}
\begin{cases}
          (x+z)(y+w) = a,\\[.5pc]
          2xz-y^2-w^2 = b,\\[.5pc]
          (x+z)(w-y) = c,\\[.5pc]
          x^2 - z^2 = d
\end{cases}
\label{eq:ilgob2}
\end{equation}
has at least one solution
$(x(a,b,c,d), y(a,b,c,d), z(a,b,c,d), w(a,b,c,d))$
for any given $a, b, c, d \in \R$.

It follows from (\ref{eq:ilgob2}) that $x^2 + y^2 + z^2 + w^2 =
\displaystyle\sqrt{a^2 + b^2 + c^2 + d^2}$. According to (\ref{eq:yosot2}) and
(\ref{eq:ilgob2}), we obtain
\begin{equation}
f(a,b,c,d) = \sigma_1(x, y, z, w)
             f\left( \sqrt{a^2 + b^2 + c^2 + d^2}, 0, 0, 0 \right)
\label{eq:yodolb2}
\end{equation}
for any $a, b, c, d \in \R$, where we denote the solution
of (\ref{eq:ilgob2}) by $(x,y,z,w)$.
Taking (\ref{eq:yodolb2}) into  account, we may introduce another
function $\sigma : \R^4 \to \{ \pm 1 \}$ that satisfies the equality
\begin{equation}
f(a,b,c,d) = \sigma(a,b,c,d) f\left( \sqrt{a^2 + b^2 + c^2 + d^2}, 0, 0, 0 \right)
\label{eq:ahob2}
\end{equation}
for all $a, b, c, d \in \R$.

Finally, define a function $m : \R \to \R$ by $m(x) = f(x, 0, 0, 0)$ for
every $x \in \R$. Then, (\ref{eq:dool2}) and (\ref{eq:ahob2}) ensure
that $m$ is a multiplicative function and that
\begin{equation*}
f(x, y, z, w) =
\sigma(x, y, z, w) \, m \left( \sqrt{x^2 + y^2 + z^2 + w^2\,}\,\right)
\end{equation*}
for all real numbers $x, y, z, w$.
\end{proof}

We will now study a stability problem of the functional equation
(\ref{eq:hana2}).
In view of Theorem~4, we can guess that the stability problem
of (\ref{eq:hana2}) is strongly connected with multiplicative
functions.

\begin{theo}[\!]
Let $X$ be a field and $K_i, L_i, M_i, N_i : X \to
[0, \infty)$ be functions for $i = 1, 2$. If a function $f : X^4 \to \C$
satisfies the following inequality
\begin{align}
&| f(x_1, y_1, z_1, w_1) f(x_2, y_2, z_2, w_2) \nonumber\\[.2pc]
&\ \ - f(x_1 x_2 + y_1 y_2 + z_1 z_2 + w_1 w_2,
x_1 y_2 - y_1 x_2 + z_1 w_2 - w_1 z_2,\nonumber\\[.2pc]
&\qquad\ x_1 z_2 - y_1 w_2 - z_1 x_2 + w_1 y_2,\,
x_1 w_2 + y_1 z_2 - z_1 y_2 - w_1 x_2) |\nonumber\\[.2pc]
   &\ \  \leq  \min\{ K_1(x_1), K_2(x_2), L_1(y_1), L_2(y_2),
                         M_1(z_1), M_2(z_2), N_1(w_1), N_2(w_2) \}\label{eq:nabi1}
\end{align}
for all $x_i, y_i, z_i, w_i \in X,$ then $f(x, 0, 0, 0)$ is either
bounded or multiplicative. Further it satisfies
\begin{align*}
&| f(x,y,z,w)^2 - f(x^2+y^2+z^2+w^2,0,0,0) |\\[.5pc]
   &\quad \ \leq
            \min\{ K_1(x), K_2(x), L_1(y), L_2(y), M_1(z), M_2(z),
                   N_1(w), N_2(w) \}
\end{align*}
for any $x,y,z,w \in X$.
\end{theo}

\begin{proof}
With $y_1 = y_2 = z_1 = z_2 = w_1 = w_2 = 0$,
(\ref{eq:nabi1}) implies
\begin{align*}
&| f(x_1, 0, 0, 0) f(x_2, 0, 0, 0) - f(x_1 x_2, 0, 0, 0) |\\[.5pc]
   &\quad \ \leq  \min\{ K_1(x_1), K_2(x_2), L_1(0), L_2(0),
                         M_1(0), M_2(0), N_1(0), N_2(0) \}
\end{align*}
for $x_1, x_2 \in X$. If we substitute $m(x)$ for $f(x,0,0,0)$ in the
above inequality, we have
\begin{align*}
&| m(x_1) m(x_2) - m(x_1 x_2) |\\[.5pc]
   & \quad \ \leq
            \min\{ K_1(x_1), K_2(x_2), L_1(0), L_2(0),
                   M_1(0), M_2(0), N_1(0), N_2(0) \}
\end{align*}
for all $x_1, x_2 \in X$.

Applying a theorem of Sz\'{e}kelyhidi \cite{13} (see Corollary 8.4
in \cite{8}), we conclude that $m$ is either bounded or multiplicative.

Finally, put $x_1 = x_2 = x$, $y_1 = y_2 = y$, $z_1 = z_2 = z$ and $w_1
= w_2 = w$ in (\ref{eq:nabi1}) to get
\begin{align*}
&| f(x,y,z,w)^2 - f(x^2+y^2+z^2+w^2,0,0,0) |\\[.5pc]
   &\quad\ \leq
            \min\{ K_1(x), K_2(x), L_1(y), L_2(y), M_1(z), M_2(z),
                   N_1(w), N_2(w) \}
\end{align*}
for all $x, y, z, w \in X$.
\end{proof}

\section*{Acknowledgement}

The first author was supported by Korea Research Foundation Grant,
KRF-DP0031.

\end{document}